\documentclass[a4paper, reqno,12pt]{amsart}
\usepackage[latin1]{inputenc}
\usepackage{latexsym}
\usepackage{amsmath,amsthm,amssymb,amsfonts}
\usepackage{enumerate}
\usepackage{graphicx}
\usepackage{hyperref}

\usepackage{graphicx}

\theoremstyle{definition}
\newtheorem{definition}{Definition}[section]

\newtheorem{example}[definition]{Example}

\theoremstyle{remark}
\newtheorem{note}[definition]{Remark}

\theoremstyle{plain}
\newtheorem{lemma}[definition]{Lemma}

\newtheorem{prop}[definition]{Proposition}

\newtheorem{theorem}[definition]{Theorem}
\newtheorem{thm}[definition]{Theorem}
\newtheorem{corr}[definition]{Corollary}

\newcommand{\M}{\mathcal{M}}

\def\setsuchas#1#2{\left\{\,{#1}\,\vrule\,{#2}\,\right\}}
\newcommand{\set}[1]{{\{#1\}}}
\newcommand{\vektor}[1]{{\boldsymbol{#1}}}
\newcommand{\Nat}{{\mathbb{N}}}
\newcommand{\Z}{{\mathbb{Z}}}
\newcommand{\C}{{\mathbb{C}}}

\newcommand{\Primes}{{\mathbb{P}}}
\newcommand{\PP}{{\mathbb{PP}}}

\newcommand{\UNI}{\mathcal{A}}
\newcommand{\divides}[2]{{#1 \left\lvert {#2} \right.}}

\newcommand{\ddivides}[2]{{#1 \left\lvert\lvert {#2} \right.}}
\newcommand{\sqbij}{\Phi}

\newcommand{\GTRUNC}{\mathfrak{S}}

\newcommand{\norm}[1]{{ \left\lvert {#1} \right\rvert}}
\newcommand{\tdeg}[1]{\norm{#1}}

\newcommand{\supp}{\mathrm{supp}}

\newcommand{\Tor}{\mathrm{Tor}}
\newcommand{\elem}{\vektor{u}}
\newcommand{\selem}{\vektor{v}}
\newcommand{\hompol}{\mathfrak{P}}

\begin{document}

\title[The ${[n]}$-truncation]{The  ring of arithmetical 
  functions with unitary convolution: the \([n]\)-truncation}
\author{Jan Snellman}
\address{Department of Mathematics\\
Stockholm University\\
SE-10691 Stockholm,
Sweden}

\email{jans@matematik.su.se}


\date{August 23, 2002}

\maketitle

\begin{abstract}
  We study a certain truncation \(\UNI_{[n]}\) of the ring of
  arithmetical functions 
  with unitary convolution,  consisting of functions vanishing on
  arguments \(>n\).  
  The truncations \(\UNI_{[n]}\) are  artinian
  monomial quotients of a polynomial ring in finitely many
  indeterminates, and are  isomorphic to the ``artinified''
  Stanley-Reisner ring \(\C[\overline{\Delta([n])}]\) of a
  simplicial complex \(\Delta([n])\).
\end{abstract}

\begin{section}{Introduction}
  The set \(\UNI\) of all complex sequences \((c_i)_{i=1}^\infty\) is in a
  natural way a complex vector space. There are several ways an
  algebra structure can be introduced on this vector space: the most
  intuitive way are through the so-called \emph{regular
    convolutions} of Narkiewicz \cite{Nark:conv}. Among those, the
  \emph{Dirichlet convolution} and the \emph{unitary convolution}
  are the best known specimen. They are in a way the two extremes, in
  that they form the maximal and minimal elements in a certain natural
  partial order on the the set of all regular convolutions.

  The ring given by Dirichlet convolution is isomorphic to the formal power
  series ring on countably many indeterminates, hence a domain; in
  fact, it is a UFD, as proved by Cashwell and Everett
  \cite{NumThe}. As a contrast,  the ring given by Dirichlet
  convolution is an epimorphic image of the formal power
  series ring on countably many indeterminates: it has zero-divisors
  and nilpotent element. In \cite{Snellman:UniDivTop} we gave a
  conjectural characterization of the zero-divisors. We also
  established some divisorial properties: a given element can have
  factorizations (into irreducibles) of different length, but there is
  always a bound for those lengths.

  Let us denote by  \(\UNI_{[n]}\) the subset of \(\UNI\) consisting
  of those sequences \((c_i)_{i=1}^\infty\)  for which \(c_i=0\) when
  \(i>n\). Restricting any convolution product on \(\UNI\), and
  modifying it so that the product of two elements in \(\UNI_{[n]}\)
  stays in \(\UNI_{[n]}\), we get a new algebra, which is an artinian monomial
  quotient of the original one. These quotients, for the case of the
  Dirichlet convolution, were studied in \cite{Snellman:Dirichlet}. It
  turned out that in this case the defining ideals are \emph{strongly
    stable}, hence we can apply the Eliahou-Kervaire resolution
  \cite{Eliahou:MinRes} to get information about various homological
  invariants.

  In this article, we study the algebra \(\UNI_{[n]}\) equipped with
  the multiplication given by the unitary convolution. It turns out
  that the defining ideals, apart from beeing monomial ideals, are
  \emph{almost square-free}, i.e. they are the sum of a square-free
  ideal and the ideal generated by the squares of the
  indeterminates. The ring is therefore the quotient  of
  \(\frac{\C[x_1,\dots,x_r]}{(x_1^2,\dots,x_r^2)}\) by a square-free
  monomial ideal, hence it is the
\emph{artinified Stanley-Reisner ring} \(\C[\overline{\Delta}]\), in
the sense of Sk\"oldberg \cite{Skold:Golod}, of a simplicial
  complex \(\Delta\). In short, the artinifed Stanley-Reisner ring is
  like the indicator algebra \(\C\{\Delta\}\) (see for instance
  \cite{Aramova:Gotzman} for a definition) but it is a quotient of
  \(\frac{\C[x_1,\dots,x_r]}{(x_1^2,\dots,x_r^2)}\) rather than a
  quotient of the exterior algebra, so it is commutative rather than
  skew-commutative. 

  The simplicial complex \(\Delta\) is easy to describe: the
  underlying vertex set consists of all prime powers \(p^a \le n\),
  and a set \(\set{q_1,\dots, q_s}\) of such prime powers is a simplex
  in \(\Delta\) if and only if the \(q_i\)'s are relatively prime, and
  their product is \(\le n\). A figure is supposedly worth a thousand
  words, so the reader may want to look at Figure~\ref{fig:Delta10} on
  page~\pageref{fig:Delta10}, where the complex for \(n=10\) is
  displayed. Understanding this complex, particularly for large \(n\),
  is the purpose of this article. 
 In particular, we
    determine the dimension and the largest non-vanishing homology of
    \(\Delta\), and show that all homology groups are
    torsion-free, by exhibiting a (non-pure) shelling. We also give an
    asymptotic estimate of the socle dimension of \(\UNI_{[n]}\).

    In \cite{Snellman:UniGenTrunc} a study was made of the truncation
    of \(\UNI\) to some subset \(V \subset \Nat^+\) which is closed
    under taking unitary divisors. The results obtained for the
    corresponding algebras \(\UNI_V\), in particular for finite \(V\),
    will be of use for us in our study of the special case
    \(V=[n]=\set{1,\dots,n}\), so this article starts with a review of
    the pertinent definitions and results from
    \cite{Snellman:UniGenTrunc}.

    \begin{subsection}{Acknowledgement}
      I used the GAP-package \emph{Simplicial Homology} \cite{SH} by  Frank
      Heckenbach, Jean-Guilleaume Dumas, Dave Saunders, and Volkmar
      Welker to calculate the homology of \(\Delta\). Having noted
      that the homology was torsion-free, I wrote a small GAP
      \cite{GAP4} programme to check if it was lex-shellable:it was,
      and it was  then   easy to prove that. I also benefitted from
      the programme Macaulay 2 \cite{MACAULAY2} to calculate Poincar\'e-Betti
      series of \(\UNI_{[n]}\).
    \end{subsection}
\end{section}

\begin{section}{Preliminaries}
  This article is a direct continuation of \cite{Snellman:UniGenTrunc}, from
  which we recall some definitions and basic results.

  \begin{subsection}{The ring of arithmetical functions with unitary convolution}
  Let \(\Nat\) denote the non-negative integers and \(\Nat^+\) the
  positive integers. 
Let \(\Primes\) denote the set of primes, with \(p_i \in \Primes\) the
\(i\)'th prime.
Let \(\PP\) denote the set
of prime powers.  \(\UNI\) denotes the set of \emph{arithmetical functions},
  i.e. functions \(\Nat^+ \to \C\). It is a \(\C\)-vector space under
  point-wise addition and multiplication by scalars, and it has a
  natural topology given by the norm \[\tdeg{f} = \frac{1}{\min
    \supp(f)},\] where 
  \[\supp(f)=\setsuchas{k \in \Nat^+}{f(k) \neq 0}.\]
\(\UNI\) becomes an
  associative, commutative \(\C\)-algebra under \emph{unitary
    convolution}
    \begin{equation}
    \label{eq:unitaryconv}
    ( f\oplus g) (n) =  \sum_{\ddivides{d}{n}} f(d) g(n/d) = \sum_{d
      \oplus m = n} f(d) g(m),
  \end{equation}
  where the unitary multiplication for positive integers is defined by 
  \begin{equation}
    \label{eq:unm}
d \oplus m = 
\begin{cases}
  dm & \text{ if } \gcd(d,m)=1\\
  0 & \text{ otherwise}
\end{cases}
  \end{equation}
and where we write \(\ddivides{d}{n}\) (or sometimes \(d \leq_\oplus n\))
when \(d\) is a unitary divisor of
  \(n\), i.e. when \(n=d \oplus m\) for some \(m\). 
For any \(k \in \Nat^+\), \(e_k\) denotes the characteristic function
on \(\set{k}\). Then \(e_1\) is the multiplicative identity, the set
of all \(e_k\) is a Schauder basis for \(\UNI\), and 
    \begin{displaymath}
      e_a \oplus e_b = e_{a \oplus b} =
      \begin{cases}
        e_{ab} & \text{ if } \gcd(a,b)=1 \\
        0 & \text{ otherwise } 
      \end{cases}
    \end{displaymath}
so \(\setsuchas{e_k}{k \in \PP}\) generates a dense subalgebra of
\(\UNI\).  

Let \(Y=\setsuchas{y_{i,j}}{i,j \in \Nat^+}\) be a doubly infinite set
of indeterminates, and let \([Y]\) be the free abelian monoid on
\(Y\). Let \(\M \subset [Y]\) consist of those monomials in the
\(y_{i,j}\)'s that are \emph{separated}, i.e. can be written
\(y_{i_1,j_1} \cdots y_{i_r,j_r}\) with \(i_1 < i_2 < \cdots <
i_r\). Then \(\M\) can be regarded as a monoid-with-zero, with the
multiplication 
\begin{displaymath}
  m_1 \cdot m_2 = 
  \begin{cases}
    m_1m_2 & \text{ if } m_1 m_2 \in \M \\
    0 & \text{ if } m_1 m_2 \not \in \M
  \end{cases}
\end{displaymath}
and
\begin{equation}
  \label{eq:Phi}
  \begin{split}
  \Phi: \M &\to \Nat^+\\
  y_{i_1,j_1} \cdots y_{i_r,j_r} & \mapsto p_{i_1}^{j_1} \cdots  p_{i_r}^{j_r}
  \end{split}
\end{equation}
is a bijection which is a monoid-with-zero isomorphism, if \(\Nat^+\)
is regarded as a monoid-with-zero with unitary multiplication. From
this follows that 
\begin{equation}
  \label{eq:s}
  \UNI \simeq \C [[\M]] \simeq \frac{\C[[Y]]}{J}
\end{equation}
where \( \C [[\M]]\) and \(\C[[Y]]\) are the generalized power series
rings on \(\M\) and \([Y]\), respectively (so \(\C[[Y]]\) is the power
series ring on bi-infinitely many variables) and \(J\) is the smallest
closed ideal of \(\C[[Y]]\) which contains all \(y_{i,j}y_{i,k}\).
  \end{subsection}

  \begin{subsection}{General truncations}
For any \(V \subseteq \Nat^+\), \(\UNI_V \subseteq \UNI\) is the
\(\C\)-sub vector space of functions supported on \(V\). With the
modified multiplication
    \begin{equation}
    \label{eq:unitaryconvV}
    \begin{split}
    ( f\oplus g) (n) &=  \sum_{d \oplus_V m = n} f(d) g(m) \\
d \oplus_V m &= 
\begin{cases}
  dm & \text{ if } \gcd(d,m)=1 \text{ and if } dm \in V\\
  0 & \text{ otherwise}
\end{cases}
    \end{split}
  \end{equation}
it becomes a \(\C\)-algebra, but in general not a sub-algebra of
\(\UNI\); it is a sub-algebra if and only if  
\[a,b \in V \implies a \oplus b \in V \cup \set{0}.\]
If \(V\) contains all unitary divisors of its elements, then
the restriction map \(\UNI  \to \UNI_V\), which is always a vector
space epimorphism, is an algebra epimorphism. In particular, if \(n\)
is a positive integer, then the set \([n]=\set{1,2,\dots,n} \subset
\Nat^+\) has this property. If we denote the kernel of the restriction
map by \(\GTRUNC_V\), then
  \begin{theorem}\label{thm:mingens}
    The set
   \begin{equation}
     \label{eq:mingensMV}
     M_V = \setsuchas{e_k}{k \not \in V, \text{ but } d \in V \text{
         for all proper unitary divisors } d \text{ of } k}
   \end{equation}
   form a minimal
   generating set of an ideal \(I_V\) whose closure is \(\GTRUNC_V\).
  \end{theorem}
  \end{subsection}

  \begin{subsection}{Finite truncations}
Now suppose that \(V\) has this property, and is finite. 
 Put 
\begin{equation}
  \label{eq:YV}
  Y(V)=\setsuchas{y_{i,j}}{p_i^j \in V \cap \PP}
\end{equation}

From \cite{Snellman:UniGenTrunc}
we have that 
\begin{equation}
  \label{eq:isoiso}
  \UNI_V \simeq \frac{\C[Y(V)]}{A_V + B_V + C_V}
\end{equation}
where
\begin{equation}
  \label{eq:abcV}
  \begin{split}
 A_V &= \left \langle y_{i,j}^2 \lvert y_{i,j} \in
Y(V) \right \rangle\\
B_V &= \left\langle y_{i,j} y_{i,k} \lvert  y_{i,j},
y_{i,k} \in Y(V) \right\rangle \\
C_V &= \left\langle y_{i_1,j_i} \cdots y_{i_r,j_r} \lvert
 y_{i_\ell,j_\ell} \in Y(V), \,  i_1 < i_2 < \cdots i_r, \,
 p_{i_1}^{j_i}\cdots p_{i_r}^{j_r} \not \in V \right\rangle
  \end{split}
\end{equation}
We also have that 
\begin{equation}
  \label{eq:arSR}
  \UNI_V \simeq \C[\overline{\Delta(V)}]
\end{equation}
where \(\Delta(V)\) is the simplicial complex on the vertex set \(V
\cap \PP\) given by 
\begin{equation}
  \label{eq:sigmain}
  \sigma=\set{p_{i_1}^{j_1}, \dots,p_{i_r}^{j_r}} \in \Delta(V) \quad
  \iff \quad p_{i_1}^{j_1} \cdots p_{i_r}^{j_r} \in V
\end{equation}
and where \(\C[\overline{\Delta(V)}]\) is the ``Artinified
Stanley-Reisner ring'' \cite{Skold:Golod} on \(\Delta(V)\): it is the
artinian commutative \(\C\)-algebra with a \(\C\)-basis
\[\setsuchas{e_\sigma}{\sigma \in \Delta(V)},\] and multiplication 
\begin{equation}
  \label{eq:ASRmult}
  e_\sigma e_\tau = 
  \begin{cases}
    e_{\sigma \cup \tau} & \text{ if } \sigma \cap \tau = \emptyset\\
   0 & \text{ if } \sigma \cap \tau \neq \emptyset
  \end{cases}
\end{equation}
Thus, as a graded \(\C\)-vector space \(\C[\overline{\Delta(V)}]\) is
isomorphic to \(\C\{\Delta(V)\}\), the indicator algebra on
\(\Delta(V)\) \cite{Aramova:Gotzman}, and as a cyclic
\(\C[Y(V)]\)-module it is isomorphic to the quotient of the ordinary
Stanley-Reisner ring \(\C[\Delta(V)]\) by the ideal generated by the
squares of all variables.
  \end{subsection}

  \begin{subsection}{The \([n]\)-truncation} 
    This article is devoted to the case \(V=[n]=\set{1,2,\dots,n}\),
    for which the ring \(\UNI_V\) and the simplical complex
    \(\Delta(V)\) have some interesting properties.
  \end{subsection}
\end{section}

\begin{section}{Some special arithmetical functions}
  We'll make use of the following special arithmetical functions: For
  a  positive integer \(n\), \(\pi(n)\) is the 
  number of primes \(\le n\), and \(\pi'(n)\) the number of prime
  powers \(\le n\). Let \(\omega(n)\) denotes the number of
  disctinct prime factors of \(n\), and let (for \(k\le 0\)) \(\pi_k(n)\)
  be the number of positive integers \(\le n\) with \(k\) distinct
  prime factors, i.e. 
  \begin{displaymath}
    \sum_{k=0}^\infty {\pi_k(n)} t^k  = \sum_{j=1}^n t^{\omega(j)}.
  \end{displaymath}

  By \(\ell(n)\) we mean the unique
  integer  such that 
  \begin{equation}\label{eq:ell}
    \prod_{i=1}^{\ell(n)} p_i \le n < \prod_{i=1}^{\ell(n)+1} p_i,
  \end{equation}
We define \(v(1)=v(2)=-1\), and for \(n \ge 3\), \(v(n)\) as the unique
  integer  such that 
      \begin{equation}
      \label{eq:vis}
         \prod_{j=2}^{v(n)+1} p_j \leq n < \prod_{j=2}^{v(n)+2} p_j
  \end{equation}
  It follows that \(v(n) = \ell(2n)-2\). The values of \(\ell(n)\) and
  \(v(n)\) for small \(n\) are tabulated below.

  \begin{center}
  \begin{tabular}{|c|c|c|c|c|c|c|c|c|c|c|c|c|c|c|}
    \hline
    n &  1 &2 & 3 & 4 & 5 & 6 & \(\cdots\) & 15 & \(\cdots\) & 30 &
    \(\cdots\) & 105 & \(\cdots\) & 210 \\ \hline
    \(\ell(n)\) & 0& 1 & 1 & 1 & 1 & 2 & \(\cdots\) & 2 & \(\cdots\) &
    3 & \(\cdots\) & 3& \(\cdots\) & 4 \\ \hline
    \(v(n)\)  &-1 &-1 & 0 & 0 & 0 & 0 & \(\cdots\) & 1 & \(\cdots\) & 1&
    \(\cdots\) & 2 & \(\cdots\) & 2\\ \hline
  \end{tabular}    
  \end{center}

\end{section}

\begin{section}{The structure of \protect{$Y([n])$}}
  We now
  determine the structure of \(Y(V)\) for  \(V=[n]\).

  \begin{definition}
    For all positive integers \(n,i\), let 
    \begin{equation}
      \label{eq:lambdadef}
      \begin{split}
\lambda_i^{[n]} & = \max  \setsuchas{j}{p_i^j \leq n} \\
\vektor{\lambda}^{[n]} &= (\lambda_1^{[n]},\lambda_2^{[n]},\dots)
      \end{split}
    \end{equation}
We may regard
    \(\vektor{\lambda}^{[n]}\) as a partition, since \(\lambda_1^{[n]} \geq
    \lambda_2^{[n]} \geq \cdots\).
    We have that 
    \begin{equation}
      \label{eq:Ylambda}
Y([n]) = \setsuchas{y_{i,j} \in
      Y}{j \leq \lambda_i^{[n]}}      
    \end{equation}
  \end{definition}

  Note that \(\max(\setsuchas{i}{y_{i,1} > 0}) = \pi'(n)\).

  \begin{example}
    If \(n=30\), then 
    \begin{displaymath}
      \vektor{\lambda}^{[30]}= (4,3,2,1,1,1,1,1,1,1),       
    \end{displaymath}
    and the variables
    \(Y(V_{[30]})\) can be visualised as in
    figure~\ref{fig:Y30}. 
    \begin{figure}[hbt]
      \begin{center}
    \setlength{\unitlength}{0.8cm}
     \begin{picture}(10,5)
      \put(0,0){\line(1,0){10}}
      \put(0,1){\line(1,0){10}}
      \put(0,2){\line(1,0){3}}
      \put(0,3){\line(1,0){2}}
      \put(0,4){\line(1,0){1}}
      \multiput(0,0)(1,0){11}{\line(0,1){1}}
      \multiput(0,1)(1,0){4}{\line(0,1){1}}
      \multiput(0,2)(1,0){3}{\line(0,1){1}}
      \multiput(0,3)(1,0){2}{\line(0,1){1}}
      
      \put(0.2,0.2){\(y_{1,1}\)}
      \put(0.2,1.2){\(y_{1,2}\)}
      \put(0.2,2.2){\(y_{1,3}\)}
      \put(0.2,3.2){\(y_{1,4}\)}
      \put(1.2,0.2){\(y_{2,1}\)}
      \put(1.2,1.2){\(y_{2,2}\)}
      \put(1.2,2.2){\(y_{2,3}\)}
      \put(2.2,0.2){\(y_{3,1}\)}
      \put(2.2,1.2){\(y_{3,2}\)}
      
      \put(3.2,0.2){\(y_{4,1}\)}
      \put(4.2,0.2){\(y_{5,1}\)}
      \put(5.2,0.2){\(y_{6,1}\)}
      \put(6.2,0.2){\(y_{7,1}\)}
      \put(7.2,0.2){\(y_{8,1}\)}
      \put(8.2,0.2){\(y_{9,1}\)}
      \put(9.2,0.2){\(y_{10,1}\)}

      \put(0.3,4.2){2}
      \put(1.3,4.2){3}
      \put(2.3,4.2){5}
      \put(3.3,4.2){7}
      \put(4.3,4.2){11}
      \put(5.3,4.2){13}
      \put(6.3,4.2){17}
      \put(7.3,4.2){19}
      \put(8.3,4.2){23}
      \put(9.3,4.2){29}
    \end{picture}
       \caption{\(Y(\vektor{\lambda}^{[30]})\)}
       \label{fig:Y30}
      \end{center}
   \end{figure}
  \end{example}

  \begin{note}
    For a fixed \(i\), as \(n \to \infty\)  we have that
    asymptotically \(\lambda_i^{[n]} \sim \frac{\log(n)}{\log(p_i)}\).
    Thus, for both \(n\) and \(i\) large we get that 
    \(\lambda_i^{[n]} \sim \frac{\log(n)}{\log(i\log(i))}\).
    This relation is illustraded in the  graph 
    figure~\ref{fig:asy}, which shows
    (the start of) \(\vektor{\lambda}^{(10^{50})}\).
  \end{note}
  \end{section}

  \begin{section}{The presentation of $\UNI_{[n]}$}
    Having determined the indeterminates occuring in
    \(Y((n])\),
    we'll consider the defining ideal \(A_{[n]} + B_{[n]} +
    C_{[n]}\) of \(\UNI_{[n]} \simeq \frac{\C[Y([n])]}{(A_{[n]} +
      B_{[n]} +  C_{[n]}}\). 
    Recall that \eqref{eq:Phi} defines a bijection \(\Phi\) between
    \(\M\) and \(\Nat^+\).
  \begin{theorem}
    Let \(V=[n]\). Then
  \begin{enumerate}
\item
 The minimal generators of \(C_{[n]}\) correspond under
  \(\sqbij\) with those \(e_k\) for which 
  \begin{itemize}
  \item \(k >n\),
  \item all  proper unitary divisors   of \(k\) are \(\leq n\),
  \item \(k=\prod_{i=1}^r p_i^{a_i}\) with \(r \leq \pi'(n)\) and \(0
    \leq a_i \leq \lambda_i^{[n]}\) for \(1 \leq i \leq r\).
  \end{itemize}
\item \(A_{[n]} + B_{[n]} + C_{[n]}\) is a strongly \(N+1\)-multi-stable   monomial ideal in
  \(\C[Y([n])]\). By this, we mean the following:
  first, let \(N\) denote the number of components \(>1\) in
  \(\vektor{\lambda}^{[n]}\). Group the variables \(y_{i,j}\) in
  column \(i\) into one group, for \(1 \leq i \leq N\), and the
  remaining variables into one last group. Order the variables
  in each group using \(\sqbij\). Then for any monomial \(m \in A_{[n]} +
  B_{[n]} + C_{[n]}\), if \(y,y'\) belong to the same group, \(\sqbij(y) <
  \sqbij(y')\), and \(\divides{y}{m}\), then 
  \(m'=\frac{y'}{y} m \in A_{[n]} +  B_{[n]} + C_{[n]}\).
 \end{enumerate}
\end{theorem}
\begin{proof}
  \begin{enumerate}
  \item 
    Follows from Theorem~\ref{thm:mingens} and \eqref{eq:isoiso} and \eqref{eq:abcV}.
  \item
    If \(m \in A_{[n]} + B_{[n]}\) then it contains two variables ``from the same
    column''. Since \(y,y'\) belong to the same group, \(m'\) must
    also  contains two variables ``from the same
    column''; thus \(m' \in A_{[n]} + B_{[n]}\).

    If \(m\) is separated, then so
    is \(m'\), and \(\sqbij(m) < \sqbij(m')\) so if in addition \(m \in
    C_{[n]}\) then \(m' \in C_{[n]}\).
  \end{enumerate}
\end{proof}

        \begin{figure}[htb]
          \begin{center}
        \caption{\(\vektor{\lambda}^{(10^{50})}\)}
          \label{fig:asy}
        \includegraphics[bb = 80 80 450 700, scale=0.4,clip]{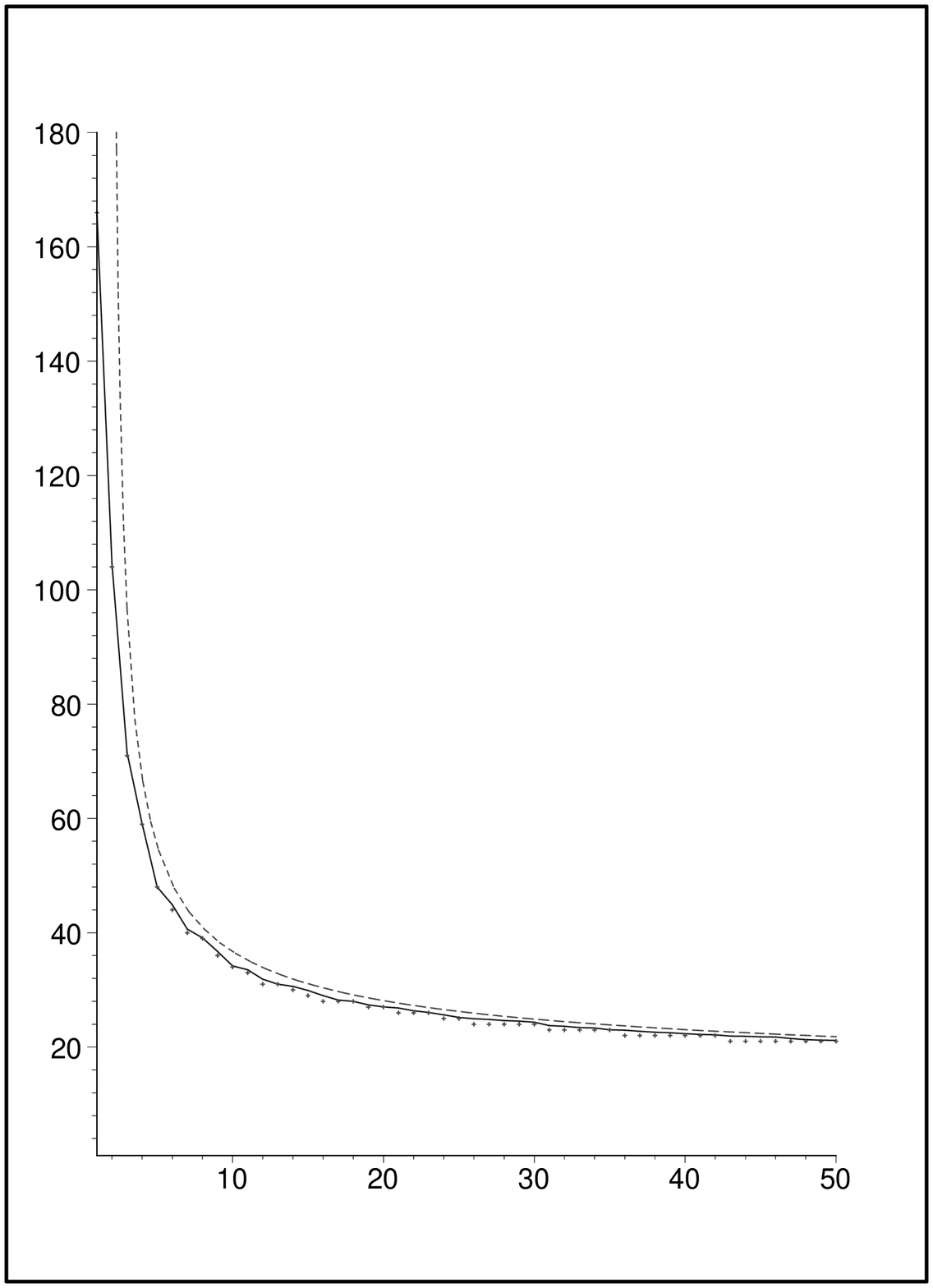}
      \end{center}
        \end{figure}

The following lemma gives a reasonably efficient way of calculating \(C_{[n]}\).
\begin{lemma}
  The minimal generators of \(C_{[n]}\)  are
contained in the set
\begin{equation}
  \label{eq:mingens}
  \M \bigcap \setsuchas{y_{i,j} \sqbij^{-1}(k)}{p_i^j \leq n, \, 1 \leq k \leq n}
  \setminus 
  \sqbij^{-1}(\set{1,2,\dots,n}) 
\end{equation}
\end{lemma}
\begin{proof}
  Take \(m=m_0 \in \M \cap Y([n])\), \(\sqbij(m)>n\). 
  Thus \(m\) is a separated monomial in variables \(y_{i,j}\) with
  \(\sqbij(y_{i,j}) \leq n\). Hence
  there is some
  \(y_{i_1,j_1}\) with \(\sqbij(y_{i,j}) \leq n\) which divides \(m\), 
  say \(m = y_{i_1,j_1} m_1\), \(m_1 \in \M \cap Y([n])\). If
  \(\sqbij(m_1) \leq n\) we are done, otherwise \(m_1\), beeing a
  separated monomials in variables \(y_{i,j}\) with
  \(\sqbij(y_{i,j}) \leq n\), can be written \(m_1 =  y_{i_1,j_2}
  m_2\), et cetera.
\end{proof}

\begin{example}
For \(n=10\), \(\vektor{\lambda}^{(10)} = (3,2,1,1)\), and
\begin{align*}
  Y(V_{[10]})) &=
\set{y_{1,1}, \, y_{1,2}, \, y_{1,3}, \, y_{2,1}, \, y_{2,2}, \,
  y_{3,1}, \, y_{4,1}} \\ 
A_{10} &= (
y_{1,1} y_{1,1}, \,
y_{1,2} y_{1,2}, \,
y_{1,3} y_{1,3}, \,
y_{2,1} y_{2,1}, \,
y_{2,2} y_{2,2}, \,
y_{3,1} y_{3,1}, \,
y_{4,1} y_{4,1}) \\
B_{10} & =(
{y_{1,1}}^2, \,
{y_{1,1}}^2, \,
{y_{1,2}}^2,\,
{y_{2,1}}^2) \\
C_{10} & = (
y_{1,1} y_{2,2}, \,
y_{1,1} y_{4,1}, \,
y_{2,2} y_{3,1}, \,
y_{1,2} y_{2,1}, \,
y_{1,3} y_{2,1}, \,
y_{2,1} y_{3,1}, \,
y_{2,1} y_{4,1}, \\ & \qquad
y_{1,2} y_{2,2}, \,
y_{1,2} y_{3,1}, \,
y_{1,2} y_{4,1}, \,
y_{1,3} y_{3,1}, \,
y_{1,3} y_{4,1}, \,
y_{2,2} y_{4,1}, \,
y_{1,3} y_{2,2}, \,
y_{3,1} y_{4,1}, \,\\ & \qquad
y_{1,1} y_{2,1} y_{3,1}, 
y_{1,1} y_{2,1} y_{4,1}, \,
y_{1,1} y_{2,2} y_{3,1}, \,
y_{1,1} y_{3,1} y_{4,1})
\end{align*}
where the last four generators are superfluous.
\end{example}

In Corollary~\ref{corr:quad}, we show that the maximal degree of a minimal
generator of \(C_{[n]}\) is \(v(n)\).

  \end{section}

  \begin{section}{Properties of $\Delta([n])$}
As an example of $\Delta([n])$, \(\Delta([10])\) is shown in
Figure~\ref{fig:Delta10}. 
\begin{figure}[htb]
      \caption{\(\Delta([10])\)}
    \begin{center}
    \setlength{\unitlength}{1truecm}
      \begin{picture}(4,4)
        \put(0,1){\circle*{0.1}}
        \put(1,1){\circle*{0.1}}
        \put(3,1){\circle*{0.1}}
        \put(4,1){\circle*{0.1}}
        \put(0,0.5){4}
        \put(1,0.5){7}
        \put(3,0.5){8}
        \put(4,0.5){9}

        \put(2,1){\circle*{0.1}}
        \put(0,3){\circle*{0.1}}
        \put(4,3){\circle*{0.1}}
        \put(2,0.5){2}
        \put(0.1,3.1){3}
        \put(3.6,3.1){5}

        \put(2,1){\line(-1,1){2}}
        \put(2,1){\line(1,1){2}}
        \put(1,2.2){6}
        \put(2.4,2.2){10}
      \end{picture}
      \label{fig:Delta10}
    \end{center}
  \end{figure}
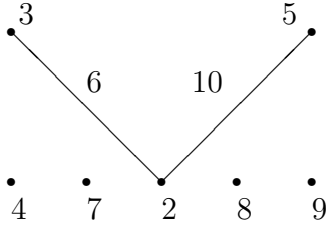


  \begin{subsection}{Simple properties}

  Let  \(r=\pi'(n)\). Clearly \(\Delta([n]\) has \(r\) vertices.
  \begin{lemma}
    \(\Delta([n])\) consists of some
  isolated vertices, and of one large component containing everything
  else. 
  \end{lemma}
  \begin{proof}
    If the vertex \(p^a\) is not isolated, but connected to \(q^b\),
    we can assume that \(p > q\). Furthermore, since \(\set{p^a,q^b}
    \in \Delta([n])\), it follows that \(p^aq^b \le n\), so \(p^aq \le
    n\), hence \(p^a\) is connected to \(q\). If \(q > 2\) then \(p >
    q > 2\) so \(\set{p^a,2} \in \Delta([n])\). Thus \(p^a\) is
    connected to 2.
  \end{proof}

  \begin{lemma}
    \(\dim \Delta([n])=\ell(n)-1\).
  \end{lemma}
\begin{proof}
   The maximal cardinality \(s\) of \((q_1,\dots,q_s)
  \in \Delta([n])\) is the maximal \(s\) such that 
  \[\prod_{i=1}^s q_i \leq n,\]
  with \(q_i \in \PP\) pair-wise relatively prime. Clearly, the best we can do is to take the
  first \(s\) prime numbers, so \(s=\ell(n)\). Hence \(\dim \Delta([n])
  = s-1 = \ell(n)-1\).
\end{proof}

\begin{lemma}
  The \(f\)-vector of \(\Delta([n])\) is
  \begin{equation}
    \label{eq:fvect1}
(f_{-1},f_0,f_1,\dots,f_{\ell(n)-1}) =
 (\pi_0(n),\pi_1(n),\pi_2(n),\dots,\pi_{\ell(n)}(n))
  \end{equation}
\end{lemma}
\begin{proof}
  We have that \(f_i\) is the number of simplicies
\(\sigma \in \Delta([n])\) of dimension \(i\), i.e. of cardinality
\(i+1\). Such a \(\sigma = \set{p_{j_1}^{a_1}, \dots,
  p_{j_{i+1}}^{a_{i+1}}}\) correspond to 
\(k=p_{j_1}^{a_1} \cdots  p_{j_{i+1}}^{a_{i+1}}\) with
  \(\omega(k)=i+1\). There are \(\pi_{i+1}(n)\) such simplices, so
\(f_i = \pi_{i+1}(n)\).
\end{proof}

\begin{lemma}\label{lemma:maxhomdeg}
    Let \(v(n)\) be defined by \eqref{eq:vis}. Then 
    \begin{enumerate}[(i)]
    \item The homological degree\footnote{The maximal \(i\) such that
        the \(i\)'th reduced homology group (with coefficients in
        \(\Z\)) of \(\Delta([n])\) is non-zero.} of \(\Delta([n])\) is \(v(n)\).
    \item The maximal homological degree
    of all \(\Delta([n])_U\), as \(U\) ranges among the non-empty
    subsets of \([n]\), is \(v(n)\).
    \end{enumerate}
\end{lemma}
\begin{proof}[Sketch of proof]
    When \(n=N=\prod_{j=2}^{s} p_j\), we have that 
  \begin{displaymath}
    \set{p_1,\dots,p_s} \not \in \Delta([n]),\quad \text{ but }
    \forall i: \set{p_1,\dots, \widehat{p_i}, \dots, p_s} \in \Delta([n])
  \end{displaymath}
  so we get \((s-1)\)-homology. As \(n\) increases and reaches
  \[2N=\prod_{j=1}^{s} p_j,\] this homology is killed off, but already
  when 
\[n=N \frac{p_{s+1}}{p_s} = p_{s+1} \prod_{j=2}^{s-1}
  p_j\]
  new \(s-1\) homology appears, since all
  \((s-1)\)-subsets of \(\set{p_1,\dots,p_{s-1},p_{s+1}}\)
  belong to \(\Delta([n])\), whereas the whole set doesn't. Filling in this
  homology, i.e. increasing \(n\) to \(p_1\cdots p_{s-1} \cdot
  p_{s+1}\), we have introduced new homology already at 
\[n=p_2\cdots p_{s-1} \cdot  p_{s+2}\] (use the fact that there
is a prime number between \(q\) and \(2q\) for all \(q\)), and so on,
until we reach \[n=\prod_{j=2}^{s+1} p_j,\] where \(s\)-homology
occurs.
\end{proof}

The asymptotic growth of \(\ell(n)\) as \(n \to \infty\) is
very slow, as the following lemma shows. Thus \(\dim
\Delta([n])=\ell(n)\) and \(v(n)=\ell(2n)-2\), the homological degree,
grows very slowly with \(n\).
\begin{lemma}\label{lemma:ellgrowth}
  There are positive real constants \(A,B,C,D\) such that for all \(n\),
  \begin{equation}
    \label{eq:asym}
    A \frac{\log(n)}{W(A \log(n))}    < \ell(n) <
   B \frac{\log(n)}{W(B\log(n))}  
  \end{equation}
where \(W(z)\) denotes the real-valued principal branch of the Lambert
W-function, defined as the root of \(W(z)\exp(W(z))=z\).
\end{lemma}
\begin{proof}
  There are positive real constants\footnote{See \cite{HW}}
  \(A_1,A_2\) such that 
  \[A_1 x <\sum_{p \leq x} \log(p) < A_2 x,\]
  hence
  \[A_1 p_n <\sum_{i=1}^n \log(p_i) <\sum_{i=1}^{n+1} \log(p_i) < A_2
  p_{n+1}.\] 
  There are constants\footnote{See \cite{HW}} \(B_1,B_2\) such that 
  \begin{displaymath}
    B_1 n \log(n) < p_n < p_{n+1} < B_2 (n+1) \log(n+1).
  \end{displaymath}
  Put \(m=\ell(n)\). Then
  \begin{displaymath}
    \sum_{i=1}^m \log(p_i) \leq \log(n) < \sum_{i=1}^{m+1} \log(p_i),
  \end{displaymath}
  so
  \begin{displaymath}
    A_1 B_1 m \log(m) < \log(n) < A_2 B_2 (m+1) \log(m+1).
  \end{displaymath}

  We claim that
  \begin{equation}\label{eq:niceeq}
    C m \log(m) = \log(n)
  \end{equation}
  has the solution
  \begin{displaymath}
    m=\frac{\log(n)}{CW(\frac{\log(n)}{C})}.
  \end{displaymath}
  From this claim, the assertion follows by monotonicity.

  Putting \(\log(n)=z\), \(\log(m)=a(z)\), \eqref{eq:niceeq} becomes
  \begin{displaymath}
    C e^{a(z)} a(z) =z,
  \end{displaymath}
  which has a solution \(a(z)=W(z/C)\). Hence
  \begin{displaymath}
    m = \exp(a(z)) = \exp(W(z/C)) = \frac{z/C}{W(z/C)} = 
    \frac{\log n}{C W(\frac{\log(n)}{C})}.
  \end{displaymath}
\end{proof}
    
  \end{subsection}

\begin{subsection}{The $h$-vector}

  Using a result by Fröberg, we have that 
  \begin{equation}
    \label{eq:hilbSR}
    \C[\Delta([n])](t) = \C[\overline{\Delta(n)}](\frac{t}{1-t}) = 
    \frac{h_0 + h_1t + h_2t^2 + \cdots +
      h_{\ell(n)}t^{\ell(n)}}{(1-t)^{\ell(n)}}  
  \end{equation}
  where \(h_0 + h_1 + \cdots +  h_{\ell(n)} \neq 1\).
  The vector \((h_0,\dots,h_{\ell(n)})\) is called the
  \emph{\(h\)-vector} of \(\Delta([n])\). It can be expressed in terms
  of the 
  \(f\)-vector as follows (see \cite{Stanley:CombCom}).

  \begin{lemma}
    The \(h\)-vector of \(\Delta([n])\) is 
    \begin{equation}\label{eq:hvect1}
        (h_0,h_1,\dots,h_{\ell(n)}) 
    \end{equation}
    where 
    \begin{equation}\label{eq:hvect2}
      \begin{split}
        h_k &= \sum_{i=0}^k (-1)^{k-i} \binom{\ell(n)-i}{k-i} f_{i-1} \\
        &= \sum_{i=0}^k (-1)^{k-i} \binom{\ell(n)-i}{k-i} c_{n,i}
      \end{split}
    \end{equation}
    with the convention that \(f_{-1}=c_{n,0}=1\).
    
    In particular,
  \begin{equation}
    \begin{split}
      h_0 &= 1\\
      h_1 &=   - \ell(n) + \pi'(n)  \\
      h_2 & = \binom{\ell(n)}{2} - (\ell(n)-1) \pi'(n) + c_{n,2}
    \end{split}
  \end{equation}
\end{lemma}

  \begin{example}
    Let \(\Delta=\Delta([10])\). Then \(\Delta\) is 2-dimensional, 
    looks like Figure~\ref{fig:Delta10}, and
    has    \(f\)-vector \((f_{-1},f_0,f_2)=(1,7,2)\)
    and \(h\)-vector \((h_0, h_1, h_2) = (1, 5, -4)\).
    Furthermore
    \[\C[\overline{\Delta}](t) = 1 +7t + 2t^2,\]
    and
    \[\C[\Delta](t) = 1+ t \frac{t}{1-t} + 2 \frac{t^2}{(1-t)^2} =
    \frac{1 + 5t - 4t^2}{(1-t)^2}.\]

  \end{example}

  \begin{subsubsection}{The \(h_2\) coefficient}
Clearly, \(h_1 > 0\) for \(n>2\).
Furthermore, we have (see \cite{ProbNum}) 
\begin{equation}
  \label{eq:pikas}
  \pi_k(x) \sim \frac{x}{\log x} \frac{(\log_2 x)^{k-1}}{(k-1)!}
\qquad x \to \infty
\end{equation}
Thus we can write (somewhat sloppily)
\begin{displaymath}
  \begin{split}
  h_2  &= \binom{\ell(n)}{2} - (\ell(n)-1) \pi'(n) + c_{n,2} \\
  & \approx \left( C \frac{\log(n)}{W(C \log(n))}\right)^2 -
  \left( C \frac{\log(n)}{W(C \log(n))} \right) \frac{n}{\log n}
  + \frac{n}{\log n} \log_2 n \\
  &= - \frac{C n}{W(C \log n)} + \frac{C^2 (\log n)^2}{W(C \log n)^2}
  + \frac{n \log 2}{(\log n)^2}
  \end{split}
\end{displaymath}
We claim that \(h_2 < 0\) for large \(n\), i.e. that
\begin{displaymath}
  n > \frac{C (\log n)^2}{W(C \log n)} + n \frac{\log 2}{C} \frac{W(C \log
    n)}{(\log n)^2}
\end{displaymath}
for large \(n\). Since \(n \gg C^2 (\log n)^2\), we need only show
that 
\begin{equation}\label{eq:toZ}
  \frac{W(C \log  n)}{(\log n)^2} \to 0 \qquad \text{ as }  n \to
  \infty
\end{equation}
However, since \(\ell(n) \to \infty\) as \(n \to \infty\),
\eqref{eq:asym} shows that 
\begin{displaymath}
  \frac{W(C \log  n)}{C\log n} \to 0 \qquad \text{ as }  n \to
  \infty,
\end{displaymath}
hence \eqref{eq:toZ} follows.
    
  \end{subsubsection}

  \begin{subsection}{Related questions}
We display below how \(h_2\) varies with \(n\). Note that
\(h_2\) has local maxima when \(\ell(n) < \ell(n+1)\), i.e. when \(n\)
is of the form
\(n=-1 + \prod_{k=1}^j p_k\).
Is the \(h_2\) coefficient negative for all \(n\)? 
\begin{center}
          \includegraphics[bb = 80 80 450 670, scale=0.5,clip]{h2.ps}
\end{center}

    One can also ask: what is the
    sign (and magnitude) of the \(h_k\) coefficient, for very large \(n\)?
  \end{subsection}
\end{subsection}

\begin{subsection}{Shellability}
  Recall that the maximal (with respect to inclusion) faces of a
  simplicial complex \(\Delta\) are called \emph{facets}, and that a simplicial complex
  is \emph{pure} if all its facets have the same dimension (which is
  then also the dimension of the simplicial complex itself). For a
  face \(\sigma \in \Delta\), we let \(\overline{\sigma} = 2^\sigma\),
  the set of all subsets of \(\sigma\), including \(\sigma\) itself,
  and the empty set.
  
  Björner and Wachs \cite{ShellNonPure1, ShellNonPure2} defines \(\Delta\) to be 
  \emph{shellable} if its facets can be arranged in linear order
  \(F_1,F_2,\dots, F_t\) in such a way that the subcomplex
  \(\cup_{i=1}^{k-1} \overline{F_i} \cap \overline{F_k}\) is pure and
  of dimension \(\dim F_k -1\) for all \(k=2,\dots,t\). They proved
  the following \cite[Lemma 2.3]{ShellNonPure1}:

  \begin{lemma}[Björner-Wachs]
    An order \(F_1,F_2,\dots,F_t\) of the facets of \(\Delta\) is a
    shelling if and only if for every \(i\) and \(k\) with \(1 \le i
    \le k\) there is a \(j\) with \(1 \le j < k\)  and an \(x \in
    F_k\) such that \(F_i \cap F_k \subseteq F_j \cap F_k = F_k
    \setminus \set{x}.\) 
  \end{lemma}

  We define the \emph{lexicographic order} on finite subsets of
  \(\PP\) in the following way: if 
    \[\sigma={a_1,\dots,a_r} \subset \PP  \quad
    \text{ with } 1 \le a_1 < a_2 < \dots < a_r ,\]
\[\tau=\set{b_1,\dots,b_s} \in \Delta([n]),  \quad \text{ with }
    1 \le b_1 < b_2 < \dots < b_s , \]
    then \(\sigma \ge_{lex} \tau\) if either \(\tau=\emptyset\) or
    \(a_1 > b_1\) or \(\sigma \setminus \set{a_1} \ge_{lex} \tau
    \setminus \set{b_1}\).
    Note that \(\ge_{lex}\) is a \emph{Boolean term-order} in the
    sense of Maclagan \cite{BoolTerm}, so 
    \begin{equation}\label{eq:bt}
      \sigma \ge_{lex} \tau \quad \implies \quad \sigma \cup \set{w}
      \ge_{lex} \tau  \cup \set{w}     \ge_{lex} \tau.
    \end{equation}

  \begin{theorem}\label{thm:lexshellable}
    \(\Delta([n])\) is shellable; the lexicographic order on the
    facets is a shelling order.
  \end{theorem}
  \begin{proof}
    Let \(F_1,\dots,F_t\) be the facets in \(\Delta([n]\) ordered
    lexicographically. 
    Using the previous lemma, we'll show that this is a shelling. So
    pick \(1 \le i < k \le t\), and let 
    \begin{displaymath}
      \begin{split}
        F_i \cap F_k &= \set{a_1,\dots,a_r} , \quad 1 \le a_1< a_2 <
        \cdots < a_r \le n\\
        F_i &= \set{a_1,\dots,a_r} \cup \set{b_1,\dots,b_t}, 
        \quad 1 \le b_1< a_b <  \cdots < b_t \le n\\
        F_k &= \set{a_1,\dots,a_r} \cup \set{c_1,\dots,c_s},
        \quad 1 \le c_1< c_2 <    \cdots < c_s \le n \\
      \end{split}
    \end{displaymath}
    We always have that \[r=\tdeg{F_i \cap F_k} < \tdeg{F_k} = r+s,\]
    since by definition, no facet is contain in another facet. If
    \[\tdeg{F_k} - \tdeg{F_i \cap F_k} = 1,\] then we are done, by
    taking \(j=i\). 

    So suppose that 
    \[\tdeg{F_k} - \tdeg{F_i \cap F_k} \ge 2.\] 
    This means that \(s \ge 2\). Since \(i<k\), \(F_i   >_{lex} F_k\).
    We distinguish two cases: some \(c_v\) is a prime power which is
    not a prime (case 1), or all \(c_v\)'s are prime (case 2).

    \textbf{Case 1:} There is some \(c_v\) which is not a prime, thus
    \(c_v = p^\delta\) with \(\delta > 1\), \(p\) a prime. Put 
    \[G=\left(F_k   \cup \set{p} \right) \setminus \set{c_v}.\]
    Then \(G >_{lex}
    F_k\), and since \(p < c_v\), \(\prod_{u \in G} u \le n\), so \(G
    \in \Delta([n])\). Clearly, \[G \cap F_k = F_k \setminus
    \set{c_v}.\] 
    Now \(G\) need not be a facet, but it is contained in one, say \(G
    \subseteq F_j\), and by \eqref{eq:bt} it follows that \(F_j
    \ge_{lex} F_k\), whence \(j < k\). Since 
    \[F_j \cap F_k \supseteq    G \cap F_k = F_k \setminus
    \set{c_v},\]
    and since, as noted, \(F_k\) can not be contained in another
    facet, we must have that 
    \[F_j \cap F_k  = F_k \setminus \set{c_v},\]
    as desired.

    \textbf{Case 2:} Since \(F_i >_{lex} F_k\) then \(b_1 < c_1 < c_2
    < \dots < c_s\).
    If all \(c_1,\dots,c_s\) are prime, then \(\gcd(b_1,c_v)=1\) for
    \(1 \le v \le s\). Hence \[G=\left(F_k \cup \set{b_1} \right)
    \setminus \set{c_1} \in \Delta([n].\]
    As before, it sufficies to
    note that \(G >_{lex} F_k\) and that 
    \(G \cap F_k = F_k \setminus  \set{c_1}\)
    to be able to conclude that there is some \(j < k\) such that
    \(G \subseteq F_j\), \[F_j \cap F_k = G \cap F_k = F_k \setminus
    \set{c_1}.\]
  \end{proof}

  It follows that all the homology groups of \(\Delta([n])\) are torsion-free.

\end{subsection}
  \end{section}

  \begin{section}{Socle degree, the Gorenstein property, and
      symmetric Hilbert function for $\UNI_{[n]}$}
    \begin{subsection}{The socle of $\UNI_{[n]}$}

      The following easy Lemma was proved in
      \cite{Snellman:UniGenTrunc}.
    \begin{lemma}\label{lemma:nSocle}
\(\mathrm{Socle}(\UNI_{[n]})\) is spanned as a \(C\)-vector space by
the set
  \begin{multline}
    \label{eq:socle}
    \setsuchas{e_k}{1 < k \leq n, \, \, e_k
      \oplus f = \mathbf{0} 
      \text{ for 
        all } f \text{ with } f(0) \neq 0} = \\
    \setsuchas{e_k}{1 < k \leq n, \,\, kp > n \text{ for all } p \in
      \Primes \text{ such that } \gcd(k,p)=1}
  \end{multline}
    \end{lemma}
    Furthermore, the \(e_k\)'s which span the socle correspond
    precisely to the facets (maximal faces) \(\sigma \in \Delta([n])\),
    \(\sigma=\set{p_{i_1}^{a_1}, \dots , p_{i_r}^{a_r}}\),
    \(k=p_{i_1}^{a_1} \cdots  p_{i_r}^{a_r}\).

  \begin{theorem} \label{thm:socle}
    Let \(\dim_\C \mathrm{Socle}(\UNI_{[n]})\) denote the vector space
    dimension of the socle of \(\UNI_{[n]}\). Then
    \begin{equation}
      \label{eq:zeta}
      \lim_{n \to \infty} \frac{\dim_\C \mathrm{Socle}(\UNI_{[n]})}{n} =
 1 - \frac{1}{2} + \sum_{i=1}^\infty \frac{ \frac{1}{p_i} -
    \frac{1}{p_{i+1}} }{\prod_{j=1}^i p_j} 
  \approx  0.60771435951661818
    \end{equation}
  \end{theorem}
  \begin{proof}
    Let \(\Nat^{++} = \setsuchas{k \in \Nat}{k > 1}\). For all \(n,k
    \in \Nat^+\), put 
    \begin{displaymath}
      \begin{split}
I_{n,0} &= \Nat^{++} \bigcap \, \left[\frac{n}{2}, n\right]\\
I_{n,k} &= \Nat^{++} \bigcap \, \left[\frac{n}{p_{k+1}}, \frac{n}{p_{k}}\right)        
      \end{split}
    \end{displaymath}
    By  \eqref{eq:socle}, the integers  \(v \in I_{n,k}\)
    correspond to \(e_v \in \mathrm{Socle}(\UNI_{[n]})\) whenever \(pv
    >n\) for all prime numbers relatively prime to \(v\).
    When  \(2p_k > n\), \( I_{n,k} = \emptyset\).
    For the remaining \(k\)'s, we have that 
    \(I_{n,k}\) contains approximatively 
    \[\frac{n}{p_k} - \frac{n}{p_{k+1}}\] integers. 
    Of those integers, only those that
    are divisible by \(p_1,\dots,p_k\) correspond to \(e_k \in
    \mathrm{Socle}(\UNI_{[n]})\). Thus, the
    contribution to the socle from \(I_{n,k}\) is approximatively
    \[\frac{ \frac{n}{p_i} - 
    \frac{n}{p_{i+1}} }{\prod_{j=1}^i p_j}. \] 
  Clearly, all integers
  in \(I_{n,0}\) are in the socle, which gives a contribution of
  approximatively \(\frac{n}{2}\). 
  Furthermore, for the intervall \(I_{n,k}\) to contain any integers,
  it must have length \(\geq 1\), i.e. \(\frac{n}{p_i} - 
    \frac{n}{p_{i+1}} \geq 1\). In particular, we must have that
      \(\frac{n}{p_i} > 1\), that is \(p_i < n\). Hence, we need only
      consider \(\pi(n)\) such intervals.
     For each interval \(I_{n,k}\) that does
      contain integers, the error 
      \begin{displaymath}
        -1 < \left(
        \sum_{\substack{x \in I_{n,k} \cap \Nat\\ \divides{p_1 \cdots
              p_k}{x}}} 
        1
        \right)
        -
        \frac{\frac{n}{p_k} - \frac{n}{p_{k+1}}}{p_1 \cdots p_k} < 1.
      \end{displaymath}

Thus, we get that 
  \begin{equation}
    \dim_\C \mathrm{Socle}(\UNI_{[n]}) \approx \frac{n}{2} +
  \sum_{\substack{k \geq 1\\p_k \leq n}} \frac{ \frac{n}{p_k} - 
    \frac{n}{p_{k+1}} }{\prod_{j=1}^k p_j}
  \end{equation}
  with an error \(<\pi(n) \approx n/\log(n)\), 
  from which \eqref{eq:zeta} follows.
  \end{proof}

  In \cite{Snellman:UniGenTrunc} we defined the \emph{multiplicative
    syzygies} of \(\UNI_V\) as the kernel \(K_2(V)\) of the
  \(\C\)-linear map 
\begin{equation}
  \label{eq:monsyz2}
  \begin{split}
  \UNI_V^+ \otimes \UNI_V^+   & \to \UNI_V^+ \\
  f \otimes g & \mapsto fg
  \end{split}
\end{equation}
We call elements in \(K_2(V)\) of the form \(e_a \otimes e_b\)
\emph{monomial multiplicative syzygies}. 

  \begin{lemma}\label{lemma:mszAn}
    The monomial multiplicative syzygies of \(\UNI_{[n]}^+\) correspond
    to the lattice points 
    \begin{equation}
      \label{eq:mozs}
      M([n])=\setsuchas{(i,j)}{1 < i,j \le n, \, ij > n} \cup
      \setsuchas{(i,j)}{1 < i,j \le n, \, \gcd(i,j) > 1} 
    \end{equation}
    Socle elements correspond to an integer on the
    \(x\)-axis such that the column supported on it is contained in
    \(M([n])\).

    Almost all syzygies are monomial, in the sense that
    \begin{equation}
      \label{eq:soclemoninf}
      \lim_{n \to \infty} \frac{\tdeg{M([n])}}{\dim_\C K_2([n])} =1
    \end{equation}
    where \(K_2([n])\) is defined as in \eqref{eq:monsyz2}.
  \end{lemma}
  \begin{proof}
    The first two assertions are obvious. By elementary linear algebra
    we have that \[\dim_\C K_2([n]) = (n-1)^2 - (n-1).\] On the other
    hand, 
    \[\tdeg{M([n])} \ge (n-1)^2 - \int_2^n \frac{n dt}{t} \ge (n-1)^2
    - n\log n,\]
    so \eqref{eq:soclemoninf} follows.
  \end{proof}

Below we have plotted the monomial
  syzygies of \(\UNI_{[30]}\). One can see that 12 is
  in the socle.

    \begin{center}
  \includegraphics[bb= 85 180 510 600, scale=0.5]{30mz.eps}
    \end{center}
    \end{subsection}

    \begin{subsection}{The Gorenstein property for $\UNI_{[n]}$}
A graded Artinian  algebra is Gorenstein if and only if the
    socle is 1-dimensional. A direct computation shows that
    \(\UNI_{[2]} \simeq \C[t]/(t^2)\) is Gorenstein.
    We note that   for \(n>2\), 
    \(e_{n-1}\) and \(e_n\) must both belong to the socle, which is
    then at least 2-dimensional, so then \(\UNI_{[n]}\) is not Gorenstein.
    \end{subsection}

    \begin{subsection}{Symmetric Hilbert function}
  A Gorenstein Artinian algebra has a symmetric Hilbert function,
  hence \(\UNI_2(t) = 1+ t\) is symmetric. 
  Can \(\UNI_{[n]}(t)\) be symmetric for other values of \(n\)?
  If \(\UNI_{[n]}(t)\) symmetric, then \(c_{n, \ell(n)} = c_{n,0} = 1\),
  which can only 
   occur when  
   \begin{equation}\label{eqn:r}
     \prod_{i=1}^r p_i  \leq n < p_{r+1}\prod_{i=1}^{r-1} p_i
   \end{equation}
   for some \(r\): in this case, \(\ell(n)=r\) and the only integer in
   the interval \([1,n]\) which is the product of \(r\) primes is
   \(\prod_{i=1}^r p_i\). Checking these intervals for \(1 \leq r \leq 10\), we
   get the matches displayed below.

   \begin{center}
   \begin{tabular}{|c|c|c|}
     \hline
     r & n & \(\UNI_{[n]}(t)\) \\ \hline
     1  & 2 & \(1+t\) \\    
     2  & 6 & \(1 + 4t + t^2\) \\
     2  & 7 & \(1 + 5t + t^2\) \\
     2  & 8 & \(1 + 6t + t^2\) \\
     2  & 9 & \(1 + 7t + t^2\)  \\
     3  & 40& \(1 + 19t + 19t^2  + t^3\) \\ \hline 
   \end{tabular}
   \end{center}

   Furthermore, we have \cite[§§ 22.11]{HW} that the average order and the
   normal order of \(\omega(n)\) is \(\log \log n\). For \eqref{eqn:r}
   we have that 
   \begin{align*}
     \log \log n &< \log \log (p_{r+1}\prod_{i=1}^{r-1} p_i) \\
     &< \log r \log(p_{r+1}) \\
     &     = \log r + \log \log p_{r+1} \\
     & < \log r + \log ((r+1) \log 2) \\
     &= \log r + \log (r+1) + \log \log 2 \ll r/2
   \end{align*}
   whenever \(r\) is sufficiently large. If  \(\UNI_{[n]}(t)\) were
   symmetric, it should be centred around \(r/2\). Hence, for
   sufficiently large \(r\), \(\UNI_{[n]}(t)\) is not symmetric. There can
   therefore be only a finite number of \(n\) such that \(\UNI_{[n]}(t)\)
   is symmetric. We conjecture that the examples tabulated above are
   in fact all such examples.
      
    \end{subsection}

\end{section}

\begin{section}{Basic homological properties}
  \begin{subsection}{\(\UNI_{[n]}\) as a cyclic \(\C[Y([n])]\)-module}
     \(\UNI_{[n]}\) is an Artinian ring and a
    zero-dimensional module 
    over \(\C[Y([n])]\), with embedding dimension \(r= \pi'(n)\), and
    homological 
    dimension \(r\). Recall \cite{HW} that \(r=\pi'(n) \approx
    \frac{n}{\log(n)}\). 
    Furthermore, for the last Betti number we have that
    \begin{displaymath}
      \begin{split}
      \beta_r(\C[Y([n])],\UNI_{[n]}) &= \dim_\C \mathrm{Socle}(\UNI_{[n]}) \\
      & \approx
      \frac{n}{2} + 
  \sum_{\substack{k \geq 1\\2p_k \leq n}} \frac{ \frac{n}{p_i} - 
    \frac{n}{p_{i+1}} }{\prod_{j=1}^i p_j} \\
  &\approx 0.60771435951661818 n,
      \end{split}
    \end{displaymath}
    by \cite[Theorem 12.4]{Stanley:CombCom} and
    Theorem~\ref{thm:socle}.
    In fact, the bijection \[\Tor_r^{\C[Y([n])]}(\UNI_{[n]},\C) \simeq
    \mathrm{Socle}(\UNI_{[n]})\] 
    is degree-preserving, so 
    \[\beta_{r,j}(\C[Y([n])],\UNI_{[n]}) = \dim_\C \mathrm{Socle}(\UNI_{[n]})_j.\]
    
For the first betti number we have
    \begin{equation}
      \label{eq:mu}
      \beta_1(\C[Y([n])],\UNI_{[n]}) = \mu(A_{[n]} + B_{[n]} + C_{[n]}) = \mu(A_{[n]}) +
      \mu(B_{[n]}) + \mu(C_{[n]}), 
    \end{equation}
    the minimal number of generators of the defining ideal.
    Clearly
    \begin{equation}\label{eq:muAB}
      \mu(A_{[n]})=r, \qquad \mu(B_{[n]}) = \sum_{i=1}^r \binom{\lambda_i^{[n]}}{2} 
    \end{equation}

  \end{subsection}

  \begin{subsection}{\(\UNI_{[n]}\) as a cyclic \(\C[\overline{Y([n])}]\)-module}
We can also consider \(\UNI_{[n]}\) as a zero dimensional module over
\(\C[\overline{Y([n])}]\), with 
embedding dimension \(r\), and infinite homologial dimension. 
We have that 
\begin{displaymath}
  \beta_1(\C[\overline{Y([n])}],\UNI_{[n]}) = \mu(B_{[n]} + C_{[n]}) = \mu(B_{[n]}) +
  \mu(C_{[n]}). 
\end{displaymath}

\begin{lemma}
  \begin{equation}
    \label{eq:UNIPO}
    \begin{split}
      P^{\C[\Delta([n])]}_{\C[Y([n])]}(t,\vektor{u}) &=  
      \text{ the square-free part of }
      P^{\C[\overline{\Delta([n])}]}_{\C[Y([n])]}(t,\vektor{u}) \\
      &=   
      t^{-1}
      \sum_{ U \subset V} \elem_U 
      t^{\tdeg{U}} \hompol_U \\ 
    P^{\UNI_{[n]}}_{\C[\overline{Y([n])}]} (t,\vektor{u}) & =
    t^{-1} \sum_{ U \subset V} \selem_U  t^{\tdeg{U}} \hompol_U 
    \end{split}
  \end{equation}
where 
\begin{displaymath}
  \begin{split}
  \hompol_U &= \sum_{i=0}^{\tdeg{U}} t^{-i}
  \widetilde{H}^i(\Delta([n])_U,\C)   \\
  \elem_U &= \elem_U(u_1,\dots,u_r) = \prod_{j \in U}   u_j, \\
  \selem_U &= \selem_U(u_1,\dots,u_r) = \prod_{j \in U}
  \frac{u_j}{1-tu_j} 
  \end{split}
\end{displaymath}
\end{lemma}
\begin{proof}
  It follows from the work of  Gasharov, Peeva, and Welker
  \cite{lcmlattice} that 
  \begin{displaymath}
  \beta_{i,\vektor{a}}(\C[Y([n])],C[\overline{\Delta([n])}])
  = \dim_\C \widetilde{H}^{i-2}(1,\vektor{y}^{\vektor{a}})_L
  \end{displaymath}
  where \((1,\vektor{y}^{\vektor{a}})_L\) is the  order complex of the
  interval 
  \((1,\vektor{y}^{\vektor{a}})\) in the sublattice \(L\) of \(Y^*\)
  generated by the minimal generators of \(A_{[n]} + B_{[n]} + C_{[n]}\).
  Similarly, 
  \begin{displaymath}
   \beta_{i, \vektor{a}}(\C[Y([n])],C[\Delta([n])])
  = \dim_\C \widetilde{H}^{i-2}(1,\vektor{y}^{\vektor{a}})_{L'}
  \end{displaymath}
  where \((1,\vektor{y}^{\vektor{a}})_{L'}\) is the  order complex of
  the interval 
  \((1,\vektor{y}^{\vektor{a}})\) in the sublattice \(L'\) of \(Y^*\)
  generated by the minimal generators of \(A_{[n]} + B_{[n]} + C_{[n]}\).
  Since the
  square-free monomials form a sublattice \(S\) of \(Y^*\), it follows
  that \(L' = L \cap S\), hence
  \begin{displaymath}
    \beta_{i, \vektor{a}}(\C[Y([n])],C[\Delta([n])]) = 
    \begin{cases}
      0 &  \vektor{a} \text{ not square-free}\\
      \beta_{i,\vektor{a}}(\C[Y([n])],\C[\overline{\Delta([n])}]) &
      \vektor{a} \text{ square-free}
    \end{cases}
  \end{displaymath}
  The remaining results are immediate from the formulaes in the appendix.
\end{proof}

\begin{thm}\label{thm:castelnuovo}
  The Castelnuovo-Mumford regularity of the
  \({\C[\overline{Y([n])}]}\)-module \({\UNI_{[n]}}\) (or, 
  equivalently, the Castelnuovo-Mumford regularity of  the
  \({\C[Y([n])]}\)-module \({\C[\Delta([n])]}\)) 
 is  \(1+v(n)\), where \(v(n)\) is as defined in \eqref{eq:vis}.
\end{thm}
\begin{proof}
  It follows from \eqref{eq:UNIPO} and \eqref{eq:vis} that the
  Castelnuovo-Mumford regularity is 
  equal to 
  \begin{displaymath}
    1+\max \setsuchas{i}{\exists U: \widetilde{H}^i(\Delta([n])_U,\C) \neq
      0} = 1+v(n).
  \end{displaymath}
\end{proof}

  \end{subsection}

  \begin{subsection}{The Stanley-Reisner ring \({\C[\Delta([n])]}\)}
    The Stanley-Reisner \({\C[\Delta([n])]}\) is a ring of dimension \(\dim
\Delta([n]) + 1=\ell(n)\) and a 
\({\C[Y([n])]}\)-module of homological dimension \(r-1\).
From the formulaes in the appendix we get that 
\begin{equation}
  \label{eq:beta1}
  \begin{split}
    \mu(B_{[n]}) + \mu(C_{[n]}) &= \beta_1(\C[\overline{Y([n])}],{\UNI_{[n]}}) \\
    & =
  \beta_1({\C[Y([n])]},\C[\Delta([n])]) \\
  &= \sum_{U 
    \subset V} \dim_\C 
  \widetilde{H}^{\tdeg{U}-2}(\Delta([n])_U,\C) \\
  &= \sum_{i=0}^{r-2} \sum_{\tdeg{U}=i+2} \dim_\C
  \widetilde{H}^i(\Delta([n])_U,\C) 
  \end{split}
\end{equation}
  \end{subsection}

  \begin{subsection}{The Koszul property}
Since \(A_{[n]}\) and \(B_{[n]}\) are quadratic, and we know for which \(i\)
the reduced 
simplicial homology \(\widetilde{H}^i(\Delta([n])_U,\C)\) can be
non-zero, we get 

\begin{corr}\label{corr:quad}
  The maximal degree of a minimal generator of \(C_{[n]}\), for \(n >
  2\), is  \(\max \set{2,  \,v(n)}\). In particular, \(C_{[n]}\) is quadratic for  \(n < 15\). 
  Furthermore, 
  \begin{equation}
    \label{eq:muC}
    \mu(C_{[n]}) = \sum_{i=0}^{r-2} \sum_{\tdeg{U}=i+2} \dim_\C
  \widetilde{H}^i(\Delta([n])_U,\C) - \sum_{i=1}^r \binom{\lambda_i^{[n]}}{2} 
  \end{equation}
\end{corr}

\begin{prop}
  \(\UNI_{[n]}\) is Koszul if and only if \(n < 15\).
\end{prop}
\begin{proof}
  \(\UNI_{[n]}\) is a monomial algebra, and is quadratic iff \(n < 15\).
Thus by a result of Fröberg \cite{Froeberg:Poincare},  \(\UNI_{[n]}\) is
Koszul iff \(n<15\). 
\end{proof}
    
  \end{subsection}

\end{section}


\appendix

\begin{section}{Homological formulaes for Stanley-Reisner rings and
    indicator algebras}

  In this appendix, we assume that \(\Delta\) is a simplicial complex
  on the finite set \(W=\set{1,\dots,r}\). For \(U \subset W\),
  \(\Delta_U=\setsuchas{\sigma \in \Delta}{\sigma \subset U}\); it is
  a simplicial complex on \(U\). We denote by
  \(\widetilde{H}^{i}(\Delta_U; \, \C)\)  the reduced simplicial
  homology. Note that when \(\emptyset \in \Delta\),
  \(\widetilde{H}^{-1}(\Delta_\emptyset; \, \C) \simeq \C\).
  
  We write
  \(S=\C[x_1,\dots,x_r]\),
  \(\bar{S}=\frac{\C[x_1,\dots,x_r]}{(x_1^2,\dots,x_r^2)}\), 
and \(E\) for the exterior algebra on the vector space of linear forms
in \(S\). Then \(\C[\Delta]\) is
a cyclic \(S\)-module, \(\C[\overline{\Delta}]\) a cyclic
\(\bar{S}\)-module, and \(\C\{\Delta\}\) a cyclic \(E\)-module. \(C\)
is a module over all these rings.

  \begin{theorem}[Hochster, \cite{Hochster:CM}]
    Let \(\beta_i\) denote the \(i\)'th Betti number of \(\C[\Delta]\)
    (in a minimal free resolution of \(\C[\Delta]\) as an
    \(S\)-module), and let \(\beta_{i,\vektor{\alpha}}\) denote the
    corresponding multi-graded Betti number. Then
    \begin{equation}
      \label{eq:Hbett}
      \begin{split}
      \beta_{i,\vektor{\alpha}} &=  
      \begin{cases}
        0 & \text{ if } \vektor{\alpha} \text{ is not square-free},
      \\
\dim_\C     \widetilde{H}^{\tdeg{U}-i-1}(\Delta_U; \, \C) &
\text{ if } \vektor{\alpha} \text{ is  square-free with }
\supp(\vektor{\alpha})=U 
       \end{cases}
\\
      \beta_i &= \sum_{U \subset W} \dim_\C
      \widetilde{H}^{\tdeg{U}-i-1}(\Delta_U; \, \C)\\
      \end{split}
    \end{equation}
    Thus the Poincaré-Betti series is
    \begin{equation}
  \label{eq:Hochster}
  P^{\C[\Delta]}_S (t,\vektor{u}) =
    \sum_{U \subset W} \prod_{j \in U}
    u_j
    \sum_{i=-\infty}^\infty t^i  
    \dim_\C \left(
      \widetilde{H}^{\tdeg{U}-i-1}(\Delta_U; \, \C)
      \right)
\end{equation}
  \end{theorem}
Note that \(\widetilde{H}^{j}(\Delta_U; \, \C) = 0\) for \(j < -1\),
\(j \geq n-1\), so the above sum is finite.

\begin{theorem}[Aramova-Herzog-Hibi \cite{Aramova:Gotzman}]
    Let \(\beta_i\) denote the \(i\)'th Betti number of \(\C\{\Delta\}\)
    (in a minimal free resolution of \(\C\{\Delta\}\) as an
    \(E\)-module), and let \(\beta_{i,\vektor{\alpha}}\) denote the
    corresponding multi-graded Betti number. Then
    \begin{equation}
      \label{eq:Arbett}
      \beta_{i,\vektor{\alpha}} =  
\dim_\C     \widetilde{H}^{\tdeg{\vektor{\alpha}}-i-1}(\Delta_{U}; \, \C), \quad
\text{ where }  U=\supp(\vektor{\alpha})
    \end{equation}

    \begin{corr}
      In the above situation, 
      \begin{equation}
        \label{eq:Arnongrad}
      \beta_i = \sum_{U \subset W} \sum_{\ell=-1}^{\tdeg{U}-1}
      \binom{\ell + i}{\ell + 1 + i - \tdeg{U}}
      \dim_\C
      \widetilde{H}^{\tdeg{U}-i-1}(\Delta_U; \, \C)
      \end{equation}
    \end{corr}
  
\end{theorem}

  \begin{lemma}[Sköldberg  \cite{Skold:Golod}]
    \label{lemma:subsrel}
    \begin{equation}
  \label{eq:subs}
 P^{\C\{\Delta\}}_E(t, u_1,\dots,u_r) =
 P^{\C[\Delta]}_S(t, \frac{u_i}{1-tu_1}, \dots, \frac{u_r}{1-tu_r}) 
\end{equation}
  \end{lemma}

\begin{corr}
  \begin{equation}
    \label{eq:poinSim}
    P^{\C\{\Delta\}}_E (t,\vektor{u}) =
    \sum_{U \subset W} \prod_{j \in U}
    \frac{u_j}{1-t u_j}
    \sum_{i=-\infty}^\infty t^i  
    \dim_\C \left(
      \widetilde{H}^{\tdeg{U}-i-1}(\Delta_U; \, \C)
      \right)
  \end{equation}
\end{corr}

\begin{definition}
For \(U \subset W\), we introduce the notation 
\begin{equation}
  \label{eq:notsa}
  \begin{split}
    \elem_U &= \elem_U(u_1,\dots,u_r) = \prod_{j \in U}   u_j, \\
    \selem_U &= \selem_U(u_1,\dots,u_r) = \prod_{j \in U}
    \frac{u_j}{1-tu_j}, \\
    \hompol_U &= \sum_{i=-1}^\infty t^{-i} \widetilde{H}^i(\Delta_U,\C)
  \end{split}
\end{equation}
\end{definition}

\begin{corr}\label{corr:POIN}
\begin{equation}
  \label{eq:poinSi}
  \begin{split}
    P^{\C[\Delta]}_S (t,\vektor{u}) & =
    t^{-1} \sum_{ U \subset W} \elem_U  t^{\tdeg{U}} \hompol_U \\
    P^{\C\{\Delta\}}_E (t,\vektor{u}) & =
    t^{-1} \sum_{ U \subset W} \selem_U  t^{\tdeg{U}} \hompol_U \\
  \end{split}
\end{equation}
\end{corr}

  \begin{lemma}[Sköldberg, \cite{Skold:Golod}]
    \label{lemma:isograd}
    \begin{equation}
      \label{eq:isograd}
      \begin{split}
       P^{\C[\overline{\Delta}]}_{\bar{S}} (t,\vektor{u}) &= 
       P^{\C\{\Delta\}}_E (t,\vektor{u}) \\
       P^{\C}_{\C[\overline{\Delta}]} (t,\vektor{u}) &= 
       P^{\C}_{\C\{\Delta\}} (t,\vektor{u}) 
      \end{split}
    \end{equation}
  \end{lemma}

\end{section}


\raggedright
\bibliographystyle{hplain}
\bibliography{journals,articles,snellman}

\end{document}